\documentclass{amsart}
\usepackage{amssymb,amscd}

\usepackage{amssymb}

\usepackage[english]{babel}

\newtheorem{theorem}{Theorem}[section]

\newtheorem{proposition}[theorem]{Proposition}

\def\H{\mathbb {H}}
\def\R{\mathbb {R}}
\def\N{\mathbb {N}}
\def\T{\mathbb {T}}
\def\Z{\mathbb {Z}}
\def\cH{\mathcal H}
\def\cD{\mathcal D}
\def\cR{\mathcal R}
\newcommand{\vol}{\operatorname{vol}}
\def\loc{\it loc}

\begin{document}

\title{A counter example to the Bueler's conjecture.}

\author{Gilles Carron}
\address{Laboratoire de Math\'ematiques Jean Leray (UMR 6629), 
Universit\'e de Nantes,
2, rue de la Houssini\`ere, B.P.~92208, 44322 Nantes Cedex~3, France}
\email{Gilles.Carron@math.univ-nantes.fr}

\maketitle
\begin{abstract}
We give a counter example to a conjecture
of E. Bueler stating the equality between the DeRham cohomology of
complete Riemannian manifold and a weighted $L^2$
cohomology where the weight is the heat kernel.
\end{abstract}

\section{Introduction}

\subsection{Weighted $L^2$ cohomology :}
We first describe weighted $L^2$ cohomology and the Bueler's
conjecture. For more details we refer to E. Bueler's paper
(\cite{Bu} see also \cite{Y}).

Let $(M,g)$ be a complete Riemannian manifold and $h\in
C^\infty(M)$ be a smooth function, we introduce the measure $\mu$
:
$$d\mu(x)=e^{2h(x)}d\vol_g(x)$$ and the space of $L^2_\mu$
differential forms :
$$L_\mu^2(\Lambda^kT^*M)=\{\alpha\in L_{\loc}^2(\Lambda^kT^*M),
\|\alpha\|_\mu^2:=\int_M |\alpha|^2(x) d\mu(x) <\infty\}.$$ Let
$d^*_\mu=e^{-2h}d^*e^{2h}$ be the formal adjoint of the operator
$d\,:\,C^\infty_0(\Lambda^kT^*M)\rightarrow
L_\mu^2(\Lambda^{k+1}T^*M)$. The $k^{\rm th}$ space of (reduced)
$L^2_\mu$ cohomology is defined by :
$$
\H^k_\mu(M,g)=\frac{\{\alpha\in L_\mu^2(\Lambda^kT^*M),
d\alpha=0\}
}{\overline{dC^\infty_0(\Lambda^{k-1}T^*M)}}=\frac{\{\alpha\in
L_\mu^2(\Lambda^kT^*M),
d\alpha=0\}}{\overline{d\cD_\mu^{k-1}(d)}}$$ where we take the
$L^2_\mu$ closure and $\cD_\mu^{k-1}(d)$ is the domain of $d$,
that is the space of forms $\alpha\in L_\mu^2(\Lambda^{k-1}T^*M )$
such that $d\alpha\in L_\mu^2(\Lambda^kT^*M)$. Also if
$\cH_\mu^k(M)=\{\alpha\in L_\mu^2(\Lambda^kT^*M),
d\alpha=0,d^*_\mu\alpha=0\}$ then we also have $\cH_\mu^k(M)\simeq
\H^k_\mu(M)$. Moreover if the manifold is compact (without
boundary) then the celebrate Hodge-deRham theorem tells us that
these two spaces are isomorphic to $H^k(M,\R)$ the real cohomology
groups of $M$. Concerning complete Riemannian manifold, E. Bueler
had made the following interesting conjecture \cite{Bu} :

{\bf Conjecture :} {\em Assume that $(M,g)$ is a connected oriented
complete Riemannian manifold with Ricci curvature bounded from
below. And consider for $t>0$ and $x_0\in M$, the heat kernel
$\rho_t(x,x_0)$ and the heat kernel measure
$d\mu(x)=\rho_t(x,x_0)d\vol_g(x)$, then $0$ is an isolated
eigenvalue of the self adjoint operator $dd^*_\mu+d^*_\mu d$ and
for any $k$ we have an isomorphism :
$$\cH_\mu^k(M)\simeq H^k(M,\R).$$
}
E. Bueler had verified this conjecture in degree $k=0$ and according to
\cite{GW} it also hold in degree $k=\dim M$. About the topological
interpretation of some weighted $L^2$ cohomology, there is results
of Z.M .Ahmed and D. Strook and more optimal results of N. Yeganefar (\cite{AS},\cite{Y}). Here we will show that we can not hope more :
\begin{theorem} In any dimension $n$, there is a connected
oriented manifold $M$, such that for any complete Riemannian
metric on $M$ and any smooth positive measure $\mu$, the natural
map :
$$\cH_\mu^k(M)\rightarrow H^k(M,\R)$$
is not surjective for $k\not=0,n$.
\end{theorem}
Actually the example is simple take a compact surface $S$ of genus
$g\ge 2$ and $$\Gamma\simeq\Z\rightarrow \hat S\rightarrow S$$ be
a cyclic cover of $S$ and in dimension $n$, do consider
$M=\T^{n-2}\times \hat S$ the product of a $(n-2)$ torus with
$\hat S$.

\section{An technical point : the growth of harmonic forms :}
We consider here a complete Riemannian manifold $(M^n,g)$ and a
positive smooth measure $d\mu=e^{2h}d\vol_g$ on it.
\begin{proposition}\label{tempere} Let $o\in M$ be a fixed point, for $x\in M$, let
$r(x)=d(o,x)$ be the geodesic distance between $o$ and $x$, $R(x)$
be the maximum of the absolute value of sectional curvature of
planes in $T_xM$ and define $m(R)=\max_{r(x)\le R} \{|\nabla
dh|(x)+R(x)\}$. There is a constant $C_n$ depending only of the
dimension such that if $\alpha\in\cH_\mu^k(M)$ then on the ball
$r(x)\le R$ :
$$e^{h(x)}|\alpha|(x)\le C_n\frac{ e^{C_n m(2R)R^2}
}{\sqrt{\vol(B(o,2R))}} \|\alpha\|_{\mu}.$$
\end{proposition}
{\it Proof.--} If we let $\theta(x)= e^{h(x)}\alpha(x)$ then $\theta$
satisfies the equation :
$$(dd^*+d^*d)\theta+|dh|^2 \theta +2\nabla dh(\theta)-(\Delta
h)\theta=0.$$ where the Hessian of $h$ acts on $k$ forms by :
$$\nabla dh(\theta)=\sum_{i,j} \theta_j\wedge \nabla
dh(e_i,e_j){\rm int}_{e_j}\theta,$$ where $\{e_i\}_i$ is a local
orthonormal frame and $\{\theta_i\}_i$ is the dual frame. If
$\cR$ is the curvature operator of $(M,g)$, the
Bochner-Weitzenb\"ock formula tells us  that
$(dd^*+d^*d)\theta=\nabla^*\nabla\theta+\cR(\theta)$. Hence, the
function $u(x)=|\theta|(x)$ satisfies (in the distribution
sense)the subharmonic estimate :
\begin{equation}
\label{subestimate} \Delta u(x)\le C_n (R(x)+|\nabla dh|(x))u(x).
\end{equation}
Now according to L. Saloff-Coste (theorem 10.4 in \cite{SC}), on $B(o,2R)=\{r(x)<2R\}$ the ball of radius $2R$
, we have a Sobolev inequality : $\forall
f\in C^\infty_0(B(o,2R))$
\begin{equation}
\label{sobolev}  \|f\|^2_{L^{\frac{2\nu}{\nu-2}}}\le C_n\frac{R^2
e^{c_n \sqrt{ k_R} R} }{\Big(\vol(B(o,2R))\Bigr)^{2/\nu}} \|df\|^2_{L^2}\end{equation}
 where $-k_R<0$
is a lower bound for the Ricci curvature on the ball $B(o,4R)$ and
$\nu=\max (3,n)$ . With (\ref{subestimate}) and (\ref{sobolev}),
the Moser iteration scheme implies that for $x\in B(o,R)$,
$$ u(x)\le C_n\frac{ e^{C_n m(2R)R^2}
}{\sqrt{\vol(B(o,2R))}} \|u\|_{L^2(B(o,2R)}.$$ From which we
easily infer the desired estimate. {\hfill $\Box$}

\section{Justification of the example and further comments}
\subsection{Justification}
Now, we consider the manifold  $M=\T^{n-2}\times \hat S$ which is
a cyclic cover of $\T^{n-2}\times S$; let $\gamma$ be a generator
of the covering group.  We assume $M$ is endowed  with a complete
Riemannian metric and a smooth measure $d\mu=e^{2h}d\vol_g$. For
every $k\in \{1,..., n-1\}$ we have a $k$-cycle $c$ such that
$\gamma^l(c)\cap c=\emptyset $ for any $l\in \Z\setminus\{0\}$
and a closed $k$-form $\psi$ with compact support such that
$\int_c \psi=1$ and such that $\Bigl({\rm support\,} \psi\Bigr) \cap \Bigl({\rm
support}\, (\gamma^l)^*\psi\Bigr)=\emptyset$ for any $l\in
\Z\setminus\{0\}$. Let $a=(a_p)_{p\in \N}$ be a non zero sequence of
real number : then the $k$-form $\psi_a=\sum_{p\in \N}
a_p(\gamma^p)^*\psi$ represents a non zero $k$ cohomology class of
$M$, indeed $\int_{\gamma^pc}\psi_a=a_p$.  We define
$R_p=max\{r(\gamma^l(x)),\ x\in c,\ l=0,...,p\}$, then if the deRham
cohomology class of $\psi_a$ contains $\alpha\in \cH^k_\mu(M)$,
then according to (\ref{tempere}), we will have
$|a_p|=\left|\int_{\gamma^pc}\alpha\right|\le M_p \|\alpha\|_{\mu} $ ;
where
$$M_p=\vol_g(\gamma^p(c))C_n\frac{ e^{C_n m(2R_p)R_p^2}
}{\sqrt{\vol(B(o,2R_p))}} \max_{r(x)\le R_p} e^{-h(x)}.$$ As a
consequence,  for the sequence defined by $a_p=(M_p+1)2^p$,
$\psi_a$ can not be represented by a element of $\cH^k_\mu(M)$.

\subsection{Further comments}
Our counter example doesn't exclude that this conjecture hold for
a complete Riemannian metric with bounded curvature, positive
injectivity radius on the interior of a compact manifold with
boundary.

\end{document}